\documentclass[11pt]{amsart}
\usepackage{geometry}                
\usepackage[parfill]{parskip}    
\usepackage{graphicx}
\usepackage{amssymb}
\usepackage{epstopdf}
\usepackage[russian, english]{babel}
\title{Did Lobachevsky have a model of his ``Imaginary geometry''?}
\author{Andrei Rodin}
\begin{document}
\maketitle
\begin{abstract}
The invention of non-Euclidean geometries is often seen through the optics of Hilbertian formal axiomatic method developed later in the 19th century. However such an anachronistic approach fails to provide a sound reading of Lobachevsky's geometrical works. Although the modern notion of model of a given theory has a counterpart in Lobachevsky's writings its role in Lobachevsky's geometrical theory turns to be very unusual. Lobachevsky doesn't consider various models of Hyperbolic geometry, as the modern reader would expect, but uses a non-standard model of Euclidean plane (as a particular surface in the Hyperbolic 3-space).  In this paper I consider this Lobachevsky's construction, and show how it can be better analyzed within an alternative non-Hilbertian foundational framework, which relates the history of geometry of the 19th century to some recent developments in the field. 
\end{abstract}    
\section{Introduction}
\paragraph {A popular story about the discovery of Non-Euclidean geometries goes like this. Since
Antiquity people looked at Euclid's Fifth Postulate (P5) with a suspicion because unlike
other Postulates and Axioms of Euclid's Elements P5 didn't seem self-obvious. For this
reason people tried to prove P5 as a theorem (on the basis of the rest of Postulates
and Axioms). Typically they tried to prove P5 by \emph{reductio ad absurdum} taking the
negation of P5 as a hypothesis and hoping to infer a contradiction from it. However the desired contradiction didn't show up. Consequences of non-P5 were unusual but not overtly contradictory. At certain point some people
including Gauss, Bolyai and Lobachevsky guessed that non-P5 opens a door into a vast
unexplored territory rather than leads to the expected dead end. Lobachevsky (or Bolyai on some accounts)
first clearly expressed this view in public. However the issue remained
highly speculative until Beltrami in his \cite{Beltrami:1868a} found some Euclidean models of Lobachevsky's
geometry, which proved that Lobachevsky's new geometry is consistent (relatively to Euclidean geometry). Finally Hilbert in his \cite{Hilbert:1899} put things in order by modernizing Euclidean axiomatic method and clarifying the logical structure of Non-Euclidean geometries.}

\paragraph {Obviously the story is oversimplified. However my task now is not to provide it with
additional details but to challenge a basic assumption, which this simplified version of
history shares with a number of better elaborated ones. This assumption concerns
the very notion of mathematical theory, which dates back to \cite{Hilbert:1899} and goes on a par with the above story.  Here is how it is described in a popular geometry textbook \cite{Greenberg:1974}:} 

\begin{quote} [P]rimitive terms, such as ``point'', ``line'' and ``plane'' are undefined and could just as
well be replaced with other terms without affecting the validity of results. ... Despite
this change in terms, the proof of all our theorems would still be valid, because
correct proofs do not depend on diagrams; they depend only on stated axioms and the
rules of logic. Thus, geometry is a purely \emph{formal} exercise in deducing certain
conclusions from certain formal premises. Mathematics makes statements of the
form ``if ... then''; it doesn't say anything about the meaning or truthfulness of the
hypotheses.
\end{quote}
\paragraph {The following passage makes it clear how the author's notion of mathematical theory has a bearing on his interpretation of history:}

\footnote{\cite{Greenberg:1974} is a geometry textbook of college level containing some historical
and philosophical material. It may be argued that it is not appropriate to take historical
and philosophical claims contained in this book too seriously and criticise them
thoroughly. I disagree. Such books written for younger students often make explicit
certain assumptions about history and philosophy of mathematics, which in more
serious studies are often taken for granted and hidden behind further details. Since my aim here is to reconsider basics rather than criticize details \cite{Greenberg:1974} serves me as a perfect reference.}

\begin{quote} 
The formalist viewpoint just stated is a radical departure from the older notion that
mathematics asserts ``absolute truths'', a notion that was destroyed once and for all
by the discovery of Non-Euclidean geometry. This discovery has had a liberating effect
on mathematics, who now feel free to invent any set of axioms they wish and deduce
conclusions from them. In fact this freedom may account for the great increase in the
scope and generality of modern mathematics.
\end{quote}

\paragraph {That one's interpretation of history of mathematics depends on one's general views on mathematics is hardly a surprise. In some sense any history of mathematics is doomed to be anachronistic: in order to study mathematics of the past one needs to fix some general ideas about mathematics at the first place. However unless one takes such general ideas dogmatically a historical study may push one to reconsider the general ideas one starts with. The principal aim of this short study is to reconsider the Hilbertian ``formalist viewpoint'' through looking back at the history of geometry of 19th century, and more specifically at Lobachevsky's works.}

\paragraph {In the following two Sections I briefly describe Lobachevsky's work and his style of reasoning. Then I  come back to the ``formalist viewpoint'' and stress some difficulties that arise when one studies Lobachevsky from this viewpoint. Finally, I propose a different viewpoint that provides a remedy. We shall see that the
question of whether or not Lobachevsky had a model of his geometry has two different answers
none of which is of yes-or-no kind. The first immediate answer is that the question is
ill-posed since Lobachevsky didn't distinguish between theories and their models in anything like the same way, in which we do this today. The second answer is subtler and more interesting. There is indeed an aspect of Lobachevsky's work relevant to our notion of model. But if one allows for the the talk of models in this context one finds something surprising: Lobachevsky didn't look for models of the geometrical theory known by his name but used a non-standard model of Euclidean plane. We shall see that this construction, which from the ``formalist viewpoint'' looks exotic and even bizarre, is crucially important for Lobachevsky's project. The alternative viewpoint suggested in the last Section of this paper provides a more natural interpretation of this Lobachevsky's construction.}

\section{Hyperbolic intuition}

\paragraph{Lobachevsky in his writings presents his main geometrical discovery in various forms and in various different contexts. 
For the following concise presentation I take Lobachevsky's STP as the principle reference. Discussing some relevant epistemological issues I shall also refer to FG and NFG. Here is a short description of Lobachevsky's geometrical works, which explains this choice. }

\footnote{The principle edition of Lobachebsky's writings is \cite{Ëîáà÷åâñêèé:1949}. Here are English titles of Lobachevsky's geometrical works listed in the chronological order and supplied with abbreviations, which I use in this paper: \\ 1823: \emph{Geometry} (G)\\
1829-30: \emph{Foundations of Geometry} (FG)\\
1835: \emph{Imaginary Geometry} (IG)\\
1835-38: \emph{New Foundations of Geometry} (NFG)\\
1836: \emph{Application of Imaginary Geometry to some Integrals} (AIG)\\
1840: \emph{Studies in Theory of Parallels} (STP)\\
1855-56: \emph{Pangeometry} (PG)\\
First publication of French version of IG: \cite{Lobachevsky:1837}\\
First publication of German version of STP:  \cite{Lobachevsky:1840}\\
English translation of STP by G. B. Halsted is printed as a supplement to  \cite{Bonola:1955}}

\paragraph{G is an early writing published only after the author's death, which contains no material related to Non-Euclidean
geometry. FG and NFG are two author's attempts to write a fundamental geometrical
treatise covering the whole of the discipline from its foundations to its special
chapters. Lobachevsky's project of rebuilding foundations of geometry developed in
these two works doesn't reduce to what became known as Lobachevskean geometry
but also includes some other new ideas which I cannot discuss here. IG and AIG have, on
the contrary, a more limited task to present of a new analytic apparatus associated with
the Lobachevskian geometry (namely, the hyperbolic trigonometry) and to demonstrate its
power through some applications. Lobachevsky introduces here this apparatus ``by hand'' reducing its geometrical
background to minimum. STP is another shortened account of the basics of
Lobachevskean geometry, which, however, is theoretically complete: it begins with
synthetics geometrical considerations and proceeds to analytic methods. STP doesn't
cover some more specific issues (like calculation of areas and volumes) treated in FG
and in NFG but unlike IG includes the foundational synthetic part. PG is the last overview
of Lobachevskean geometry written by the author; it is less systematic than STP and
fixes some minor technical problems, which Lobachevsky found in STP after its
publication. \\ This description makes it clear that STP is the best compact systematic
presentation of the topic written by Lobachevsky himself. Importantly Lobachevsky's
notion of ``Imaginary geometry'' remains stable across all of these works.}

\paragraph {STP is written in the classical Euclidean ``synthetic'' style reinforced by analytic
methods described in Section 4. As far as the logical structure of presentation
is concerned it is apparently not of Lobachevsky's major concern. Lobachevsky
presents to the reader a list of propositions without specifying which of them are
definitions, which are assumed as axioms and which are assumed as commonly known
theorems (independent of P5); among proved propositions there are
theorems known to Lobachevsky from his sources as well as theorems first proven by
Lobachevsky himself. From a historical viewpoint these features of Lobachevsky's
style are hardly surprising since all of Lobachevsky's predecessors and
contemporaries working on the ``Problem of parallels'' also followed the same
traditional line. I stress these features only in order to confront the widespread philosophical myth according to which the invention of Non-Euclidean geometry required an abrupt departure from the ``usual'' geometrical intuition.}

\paragraph {Instead of P5 Lobachevsky uses the following Axiom of Parallels (AP) known to be equivalent to P5
since Antiquity:}
\footnote{AP is also known under the name of Playfair's Axiom}

\begin{quote}(AP) Given a line and a point outside this line there is unique other line which is parallel
to the given line and passes through the given point.\end{quote}

\paragraph {Here the term ``parallels'' stands as usual for straight lines having no common points.
We'll se shortly how Lobachevsky changes this Euclidean terminology. For a
terminological convenience I shall call a given straight line \emph{secant} of another given
straight line when the two lines intersect (in a single point). Let's now make the
required construction and listen what our intuition says about it. Although the intuition says us nothing
definite as to whether AP is true or not it says several other important things:}

\paragraph {(i) Parallel lines exist (unlike round squares); moreover through a given point $P$ outside
a given straight line $l$ passes at least one parallel line $m$. The required construction can be made on the basis of Euclid's Postulates except P5. Drop a perpendicular $PS$ from $P$ to $l$ and then produce another
perpendicular $m$ to $PS$ passing through $P$ . The fact that $m$ is parallel to $l$ follows from
the theorem about an external angle of triangle (Proposition 1.16 of Euclid's \emph{Elements}, which doesn't depend on P5 (Fig.1).}

\footnote{For suppose that $m$ and $l$ intersect in $A$. Then the external angle $RPA$ is equal to the
internal angle $PSA$. This contradicts the theorem about an external angle which implies
that $RPA$ must be strictly superior to $PSA$.}

\includegraphics {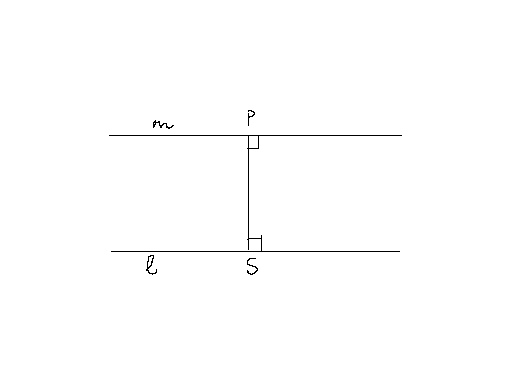}

\begin{center}
Fig.1
\end{center}

\paragraph {(ii) Given a straight line and a point outside this line there exist secants of the given
line passing through the given point. To construct a secant take any point of the given
line and connect it to the given point outside this line.}

\paragraph {(iii) Let $PS$ be perpendicular to $l$ and $A$ be a point of $l$. Consider a straight line $PR$ such
that angle $SPR$ is a proper part of angle $SPA$ (and hence is less than angle $SPA$).
Given this I shall call line $PR$ \emph{lower} than line $PA$ (and call $PA$ \emph{upper} than $PR$). Notice that
this definition involves the perpendicular $PS$, and so depends on the choice of $P$. Then
$PR$ intersects $l$ in some point $B$, i.e. it is a secant. In other words a line,
which is lower than a given secant is also a secant (Fig.2).}

\footnote{To prove (iii) rigorously one needs Pasch's axiom which Lobachevsky never mentions
but often uses tacitly. This axiom first introduced in \cite{Pasch:1882} reads: Given a triangle and a straight line intersecting one of the triangle's sides but passing through none of the triangle's apexes the given line intersects one of the two other sides of the given triangle. To apply this axiom to the given case one needs a simple additional construction that I leave to the reader. Remind that my argument here concerns common intuition but not rigorous proofs: whatever improvement on (iii) can be possibly made it remains
intuitively evident.}

\includegraphics {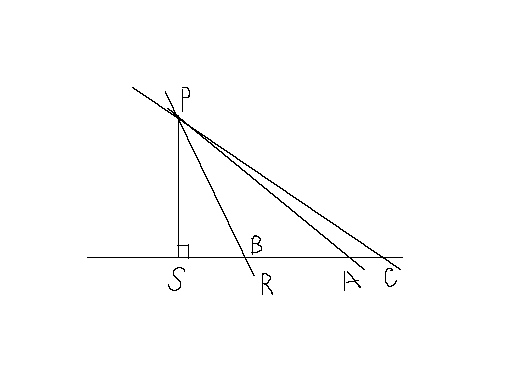}
\begin{center}
Fig.2
\end{center}

\paragraph{(iv) There exist no upper bound for secants of a given line passing through a given
point outside this given line. For given some secant $PA$ one can always take a further
point $C$ such that $A$ will lay between $S$ and $C$ and so secant $PC$ be upper than the given
secant $PA$ (Fig. 2).}

\paragraph{(v) Let $m$ be parallel to $l$ , which is constructed as in (i). Let $n$ be another parallel to $l$
passing through the same point $P$. Suppose that $n$ is lower than $m$ (obviously this
condition doesn't restrict the generality). Then any straight line which is upper than $n$
and lower than $m$ is also parallel to $l$ (Fig. 3)}

\includegraphics {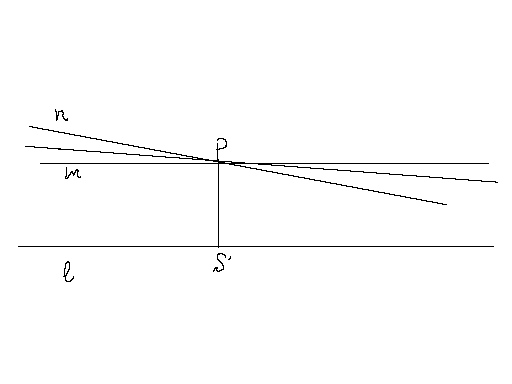}
\begin{center}
Fig.3
\end{center}

\paragraph{(vi) Parallels to a given straight line passing through a given point have a lower bound.
To assure it rigorously one needs some continuity principle like one asserting the
existence of Dedekind cuts. Then (vi) follows from (iv). Lobachevsky doesn't states
such a principle explicitly but endorses (vi) anyway.}

\paragraph{(vii) Any straight line $PA $- a secant or a parallel - passing through point $P$ as shown at
Fig. 2 is wholly characterised by its characteristic angle $SPA$. In particular this
concerns the lowest parallel mentioned in (vi). Let the measure of $SPA$ corresponding
to the case of the lowest parallel be $\alpha$. Now it is clear that by an appropriate choice of
$l$ and $P$ one can make $\alpha$ as close to $\pi/2$ as one wishes. For given any angle $SPA < \pi/2$
it is always possible to drop perpendicular $AT$ on $PS$ (Fig.4). Then $PA$ is a secant of $AT$
and so by (iii) all parallels to $AT$ including its lowest parallels are upper than $PA$. Hence
the value of $\alpha$ corresponding to straight line $AT$ and point $P$ outside this line is between
$SPA$ and $\pi/2$. Since the only variable parameter of the configuration is the distance $d$
between the given straight line and the given point outside this line the angle $\alpha$ is wholly
determined by this distance.}

\includegraphics {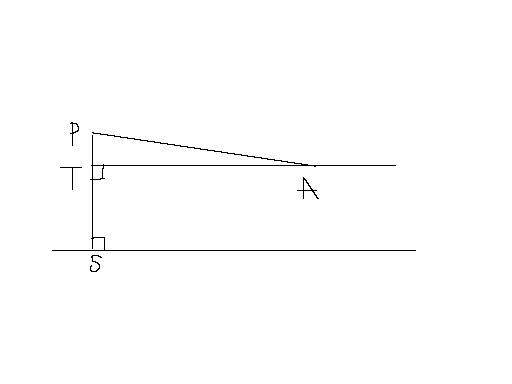}
\begin{center}
Fig.4
\end{center}

\paragraph{(i - vii) provide the intuitive basis for Lobachevsky's Imaginary geometry (see STP,
propositions 7, 16, 21). He proceeds as follows. First, he makes a terminological
change: he calls ``parallels'' (not just non-intersecting straight lines but) the two
boundary lines which separate secants from non-secants (i.e. parallels in the usual
terminology) passing through a given point. So in Lobachevsky's terms there exist
exactly two parallels to a given straight lines passing through a given point, which may
eventually coincide if AP holds (i.e. in the Euclidean case). For further references I
shall call these two parallels \emph{right} and \emph{left} (remembering that this assignment of parity is arbitrary). Since Lobachevsky's definition of parallels involves the choice of $P$ it is not immediately clear that parallels so defined form equivalence classes. So Lobachevsky must show that the property of being parallel (in his new
sense) to a given straight line is independent of this choice (STP, proposition 17), and
that the relation of being parallel is symmetric and transitive (while reflexivity may be
granted by the usual convention) (STP, propositions 18, 25). For the obvious reason
transitivity may work here only for parallels of the same orientation, i.e. separately
for right and for left parallels. Lobachevsky provides the required proofs making them in the traditional synthetic Euclidean-like manner. Then Lobachevsky proves some further properties of parallels, in particular the fact that the angle $\alpha$
characterising a parallel (see vii above) can be made not only however close to $\pi/2$ but
also however close to zero (STP, proposition 23). This immediately implies that if P5
doesn't hold then given an angle $ABC$, however small, there always exist a straight line
$l$ laying wholly inside this angle and intersecting none of its two sides (Fig.5). This is already by far more counterintuitive than (i-vii) but still not counterintuitive enough to rule out this construction as absurd and on this ground to claim a proof of P5.}

\includegraphics {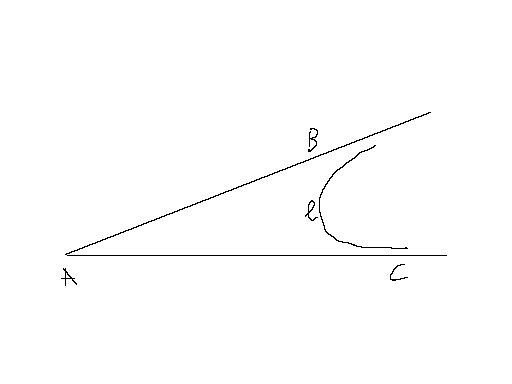}
\begin{center}
Fig.5
\end{center}

\paragraph{Lobachevsky's aim is to develop a general theory of geometry, which would not depend of P5 and include Euclidean geometry as a special case, but not to develop an alternative geometry. In the \emph{Introduction} to his NFG Lobachevsky makes this point explicitely:}

\begin{quote}
The principle conclusion, to which I arrived .... was the possibility of Geometry in a
broader sense than it has been [earlier] presented by Founder Euclid. This extended
notion of this science [=of Geometry] I called Imaginary Geometry; Usual [=Euclidean]
Geometry is included in it as a particular case.
\end{quote}

\footnote{Hereafter translations of Lobachevsky's passages from Russian are mine.}

\paragraph{Formally speaking, (i)-(vii) indeed cover the Euclidean case: in this case $\alpha = \pi/2$ and so the two Lobachevsky's two parallels coincide. What remains problematic here is the nature of variation of $\alpha$. (vi) says us nothing about the value of $\alpha$ except that it is positive but doesn't exceed $\pi/2$. Does this mean that one can stipulate by fiat any value of $\alpha$ from this interval? }

\paragraph{Geometry traditionally makes a sharp distinction between
universally valid propositions (axioms and theorems) and particular constructions with
their stipulated properties. One is free to produce a geometrical construction with any desired
properties as far as Euclid's Postulates (or some other fixed constructive principles) grant that the wanted construction is doable. For example, one is free to produce a right angle, an acute angle or an obtuse angle depending on one's personal taste or one's specific purpose. But in this traditional setting one is not free to stipulate axioms and constructive principles (postulates) in a similar way. However in Lobachevsky's reasoning this traditional distinction is blurred. According to (vii) both Lobachevsky's parallels passing through a given point $P$ are uniquely determined by the distance $d = PS$. This means that given a line $l$ and a point $P$ outside this line the angle $APS = \alpha$ has some definite value, which cannot be any longer a matter of stipulation. But since we don't know this value we can only make some hypothetical reasoning about it. If $\alpha = \pi/2$ (the Euclidean case) then the same holds for any other choice of $l$ and point $P$  (see STP, proposition 20). If $\alpha < \pi/2$ (the Hyperbolic case) the situation is more involved, because $\alpha$ depends on $d$. Here is Lobachevsky's fundamental equation, which describes this dependence:}

\begin{equation} \label{eq:1}
\tan{(\alpha/2)} = a^{-d}
\end{equation}

\paragraph{where $a$ is a positive factor. How to interpret this formula? The factor makes a new
trouble. On the one hand, it is clear that the unit used for measuring the distance $d$ can be
always chosen in such a way that \eqref{eq:1} takes this most convenient form:}
 
\begin{equation} \label{eq:2}
\tan{(\alpha/2)} = e^{-d}
\end{equation}

\paragraph{were $e$ is the base of natural logarithms. Then $\alpha$ gets determined by $d$ as expected. As
Lobachevsky says in STP this move "simplifies calculations". However such a choice of unit determines the value of $\alpha$ in each particular case! It turns out that the value of $\alpha$ depends of our arbitrary choice of
unit, which is supposed to be a matter of convention having no theoretical significance. Turning things the other way round one may also say that the unite of length is uniquely determined here by $\alpha$, that is, that the usual freedom of choice of units cannot be granted in this case. This property of hyperbolic spaces first noticed by Lambert (see \cite{Bonola:1955}) is often described as the existence of \emph{absolute unite} of length. However one shouldn't forget that this \emph{absolute unite} depends on the variable factor $a$. This actually means that unlike Euclidean geometry Lobachevskian geometry is not a single geometry but a family of similar geometries parametrized by a positive real factor. The geometrical nature of this factor remained a mystery at least until 1868 when Beltrami \cite{Beltrami:1868b} discovered the theoretical link between Lobachevsky's and Riemann's works (see Section 5 below). So it is hardly surprising that Lobachevsky's distinction between between the general theory of geometry and particular geometries covered by this general theory remains rather unclear. In the next Section we shall see how this problem was later settled by Hilbert and his followers who promoted what Greenberg calls the ``formalistic viewpoint'' in geometry.} 

 \paragraph{Although the controversy between the "formal" and the "intuitive" approaches to geometry didn't yet appear in 1830-ies when Lobachevsky made his most important discoveries a similar controversy was already known. I mean the controversy between \emph{analytic} and \emph{synthetic} ways of doing geometry. Here is what Lobachevsky says about it in the Preface to his NFG:}

\begin{quote}
In Mathematics people use two methods: analysis and synthesis. A specific
instrument of analysis are equations, which serve here as the first basis of any
judgement and which lead to all conclusions. Synthesis or the method of constructions
involves representations immediately connected in our mind with our basic concepts. ... Science starts with a pure synthesis; all the rest is produced by jugement which derives new data from the first data given by synthesis and thus broadens our knowledge unlimitedly into all directions. Without any doubt the first data are always
acquired in nature through our senses. Mind can and must reduce them to minimum, so
they could serve as a solid foundation for science.
\end{quote}

\paragraph{The fact that Lobachevsky stresses the importance of synthetic approaches in
science in general and in geometry in particular shows that unlike some of his
contemporaries he was not sympathetic to the idea of replacing intuitive
geometrical reasoning by some sort of calculus. He rather believed that the spatial
intuition and the spatial experience are ultimate sources of geometrical truths
and that analytic methods serve only for "derivation of new data from first data".}

\paragraph {Consider also this interesting passage from NFG where Lobachevsky anticipates physics of the 20th century by putting forward a view according to which geometry of physical space is variable and depends on physical factors:}

\begin{quote}
[T]he assumption according to which some natural forces follow one Geometry while
some other forces follow some other specific Geometry, which is their proper
Geometry, cannot bring any contradiction into our mind.
\end{quote}

\paragraph {This passage clearly shows that Lobachevsky thinks about alternative Non-Euclidean geometries  as prospective mathematical tools for physics but not as mere mathematical abstractions. We find in Lobachevsky no sign of Greenberg's enthusiasm about the alleged freedom of mathematicians to ``invent any set of axioms they wish and deduce conclusions from them''. Although the discovery of Non-Euclidean geometries could later motivate proponents of such a freedom Lobachevsky himself certainly didn't enjoy this freedom when he made his discovery.}   

\paragraph{Taking into consideration what has been said so far one may
come to the conclusion that the question "Did Lobachevsky have a model of his geometry?"
has no more sense than the question whether or not Euclid had a model of his geometry. In the Section 4 we shall see that in fact this question allows for a more specific and more interesting answer. But before coming to this discussion let me return to Greenberg's ``formalist viewpoint'' and show how the emergence of this viewpoint relates to Lobachevsky's work.}

\section{Hilbertian Axiomatic Method}

\paragraph {What Greenberg in the above quote describes as the ``formalistic viewpoint'' comes with a specific method of building mathematical theories, which dates back to Hilbert's treatise on foundations of geometry \cite{Hilbert:1899}    and his influential paper \cite{Hilbert:1918}. This method is often called after Hilbert simply \emph{the Axiomatic Method}. This common name hides the fact that the notion of \emph{axiom} relevant to this modern method of theory-building strikingly differs from what used to be so called before the 20th century. In order to stress the specific character of \emph{this} Axiomatic Method I call it hereafter \emph{Hilbertian}. In the rest of this Section I remind the reader basic features of Hilbertian Axiomatic Method and then show how this method solves some epistemological problems related to Non-Euclidean geometries. The purpose of this Section is not to promote the Hilbertian Axiomatic Method which hardly needs any further promotion but rather prepare a background for its critical reconsideration in the following Sections.}

\paragraph {The Hilbertian Method identifies a mathematical theory with a system of propositions some of which are assumed as \emph{axioms} while some other (called \emph{theorems}) are deduced from the axioms according to fixed rules of logical inference. Propositions (i.e. axioms and theorems) are viewed here in two different ways. First, they are viewed as syntactic constructions having no meaning and truth-value. So conceived propositions are called \emph{formal}; systems of formal propositions are called \emph{formal theories}. Formal propositions and formal theories are provided with meaning and truth-values through a special procedure of \emph{interpretation}, which assigns to terms of a given formal some particular mathematical objects. Thus formal theories and propositions become \emph{interpreted}; interpreted propositions have certain truth-values, which obviously depend on the given assignment. An assignment, which makes all provable (deducible) propositions of a given theory true is called a \emph{model} of this theory. A given theory may have multiple models and multiple "would-be-models", in which some formally provable propositions are true but some other turn to be false.}

\paragraph {The role of models is (at least) twofold. First, models provide an
intuitive support allowing, for example, for thinking of proposition ``given two points
there exist an unique straight line going through these points'' in the usual way.
(Alternatively one can think of points in the way one usually thinks of straight lines
and think of straight lines in the way one usually thinks of points. It would make a
difference in Euclidean geometry, in which there exist lines without common points, but
not in Projective geometry in which any two straight lines intersect.) Second, models
are used for proving \emph{relative consistency} of a given theory $T$ with respect to some other theory $T'$, which is called in this context a \emph{metatheory}. This works as follows. One takes $T'$ and some its model $M'$ for granted and use these things for building a model of $T$. In particular, one may take for granted Arithmetic (assuming that this theory comes with a model)  and use Arithmetic for building models of various geometrical theory. The existence of such a arithmetical model $M$ of a given geometrical theory $T$ shows that $if$ Arithmetic is consistent $then$ the theory $T$ is also consistent. If one further assumes that Arithmetic is more ``secure`` or more ``basic'' than geometrical theories (as did Hilbert) than this relative consistency result presents an important epistemic gain.     
}

\paragraph {Let me now explain the meaning of the word ``formal'' used in this context. Commenting in a letter on his \cite{Hilbert:1899} Hilbert writes:} 

\begin{quote}[E]ach and every theory can always be applied to the infinite number of systems of basic elements. One merely has to apply a univocal and \underline{reversible} one-to-one transformation [to the elements of the given system] and stipulate that the axiom for the transformed things be correspondingly similar (quoted by \cite{Frege:1971}, underlining mine)\end{quote}

\footnote {When Greenberg says that in a formalistic geometrical theory ``Primitive terms, such as ``point'', ``line'' and ``plane'' ... could just as well be replaced with other terms without affecting the validity of results'' he obviously means not only that these geometrical objects can be called differently; he also means that these geometrical objects (as traditionally conceived) can be replaced by some other objects. }   

\paragraph {Clearly Hilbert assumes here the following property of formal theories: all models of a given formal theory are transformable into each other element-wise by reversible transformations. In modern language such models are called \emph{isomorphic} while a theory having the property that all its models are isomorphic is called (after Veblen) \emph{categorical}. So the idea behind the  Hilbertean Axiomatic Method is this: an axiomatic theory captures the common \emph{form} or \emph{structure} shared by a bunch of isomorphic mathematical constructions, no matter how these constructions are produced and how they are conceived intuitively. In the current literature the term ``structure'' is used in the given context more often than the term ``form''. More precisely a structure in the relevant sense of the term can be described as a result of \emph{abstraction} from a bunch of isomorphic constructions, which amounts to ignoring all the differences between such constructions, i.e., to considering these constructions ``up to isomorphism''. The fact that ``being isomorphic'' is an equivalence relation is essential for making such an abstraction possible \cite{Rodin:2010}. The view according to which mathematics is a science of structures is known under the name of Mathematical Structuralism \cite{Hellman:forthcoming}.}     

\paragraph {We see that the Hilbertian Axiomatic Method is in fact something more than just a method. It comes with a specific epistemic framework for doing mathematics. Within this framework a (formal) mathematical theory can be described as a set of axioms taken together with its deductive closure, i.e. with all propositions deducible from the given axioms. Any such mathematical theory must be logically \emph{consistent}, i.e., not allow for deriving a contradiction. Further features of theories are desired but optional: a good theory comprises a set of axioms, which are logically independent; it is categorical; it has interesting models; it solves some earlier stated problems; it can be successfully applied in physics or elsewhere. }

\paragraph {Let us now see how the Hilbertian Axiomatic Method solves the epistemic problem concerning Lobachevskian geometry stressed in the previous Section. The key point of this solution is this: the new method of theory-building doesn't involve the traditional distinction between general geometrical principles (Postulates and Axioms) and stipulated properties of particular geometrical constructions. Even if some distinction of this sort can be still drawn within the Hilbertian framework it no longer plays any important role. In this framework one is authorized both to assume and to not assume AP just like in the traditional Euclidean framework one is authorized to assume and not to assume that a given triangle is isosceles. (As a matter of course in the Hilbertian framework one is \emph{not} authorized to assume AP and its negation simultaneously just like in the traditional framework one is not authorized to assume that one and the same triangle is and is not isosceles.) Thus the Hilbertian framework allows one to do Euclidean geometry (assuming AP) as well as Lobachevskian geometry (not assuming AP). The new framework allows for a peaceful co-existence of these geometrical theories (as different parts of the same body of mathematical knowledge) in spite of the fact that they are logically incompatible.}     

\paragraph {For the question of historical origins of Hilbertian Axiomatic Method I refer the reader to \cite{Toepell:1986} and further literature wherein. Although the discovery of Non-Euclidean geometry certainly played a role in the development of this Method this development was not a direct consequence of this discovery. Other important conceptual developments took place during the second half of the 19th century when Non-Euclidean geometries were already known but the new Axiomatic Method was not available yet. I shall say more about these developments in Section 5 below. Let us now see whether it is possible to identify some elements of the Hilbertian Method in Lobachevsky's writings.}

\section{Hyperbolic calculus}

\paragraph{In order to see that the modern notion of model is not totally irrelevant to Lobachevsky's work
consider the following quote from FG: }
\begin{quote}
The geometry on the limiting sphere is exactly the same as on the plane. Limiting
circles stand for straight lines while angles between planes of these circles stand for
angles between straight lines.
\end{quote}

\paragraph{Even without knowing the exact sense of Lobachevsky's terms (which I shall shortly
explain) one can identify here a basic element of the Hilbertean framework, namely the idea that
usual geometrical terms like ``straight line'' and ``angle between straight lines'' can be given some unusual meanings without producing any essential change in the corresponding theory (in this case - Euclidean geometry). Speaking in today's terms Lobachevsky describes here a non-standard model of Euclidean plane. Why not a model of his new Hyperbolic geometry as the today's reader would most probably expect? Why Lobachevsky translates convenient notions of Euclidean geometry into the new language of Hyperbolic geometry rather than the other way
round? Let me now explain why and how.}

\paragraph{Basic facts about Hyperbolic plane, which can be proved by traditional synthetic methods, were mostly known before Lobachevsky. Lobachevsky first managed to supply the system of synthetic reasoning briefly described above in Section 2 with an appropriate analytic tool, which allowed him to do \emph{analytic geometry} on the Hyperbolic plane.  The non-standard model of Euclidean plane serves Lobachevsky for developing this analytic apparatus. The following presentation follows STP closely.}

\paragraph{In Euclidean geometry there are two kinds of sheaves of straight lines: (a) sheaves of
parallel lines and (b) sheaves of lines passing through the same point. Given a sheaf of
either sort consider a line (or surface in 3D case) normal to each line of the given
sheaf. So you get (a) either a straight line (or plane in 3D case) or (b) a circle (or sphere in
3D case ) (Fig.6 a, b)}

\includegraphics {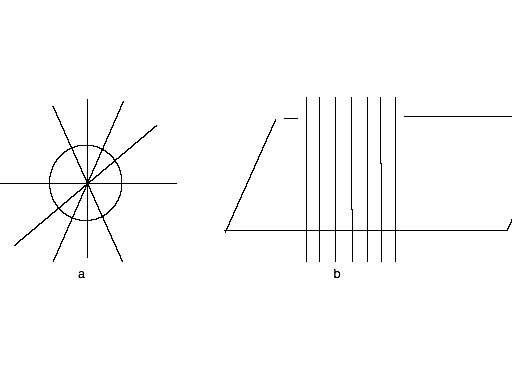}
\begin{center}
Fig.6
\end{center}

\paragraph{In Lobachevskian (Hyperbolic) geometry both configurations shown at Fig.7 a, b exist
although lines at Fig.7b are not parallels in Lobachevsky's sense. In addition one
gets here a new specific sort of sheaf, namely that of Lobachevsky's parallels.
This new sheaf comes with a new normal line and a new normal surface, which
Lobachevsky calls \emph{limiting circle} (or otherwise \emph{horocircle}) and \emph{limiting sphere}
(otherwise \emph{horosphere}) correspondingly (Fig.7):}

\includegraphics {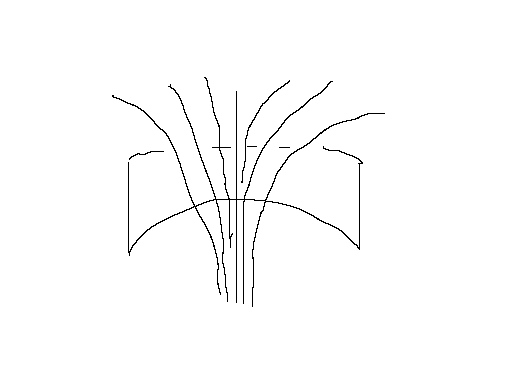}
\begin{center}
Fig.7
\end{center}

\paragraph{To see that horocircles on given horosphere verify AP (and in fact the rest of axioms
of Euclidean geometry) observe the following. Call (as usual) a given straight line $l$
parallel to a given plane $\alpha$ just in case $l$ is parallel (in Lobachevsky's sense) to its
orthogonal projection $m$ onto $\alpha$. It can be then easily shown by usual synthetic methods
(I leave it as an exercise) that for all $l$ and $\alpha$ defined as before there exist a unique plane $\beta$
having no common point with $\alpha$ (that is, parallel to $\alpha$ in the usual sense) such that $l$
lays in $\beta$. This lemma, which resembles AP in a way, doesn't depend on AP. Notice that
any horocircle laying on a given horosphere can be obtained as an intersection of the
horosphere with a plane parallel to the sheaf of Lobachevsky's parallels corresponding
to this horosphere. This immediately implies that the non-standard interpretation of
Euclidean geometry suggested by Lobachevsky verifies AP (the horocircles are called
here parallel in the usual Euclidean sense of having no shared point) (Fig. 8):}

\includegraphics {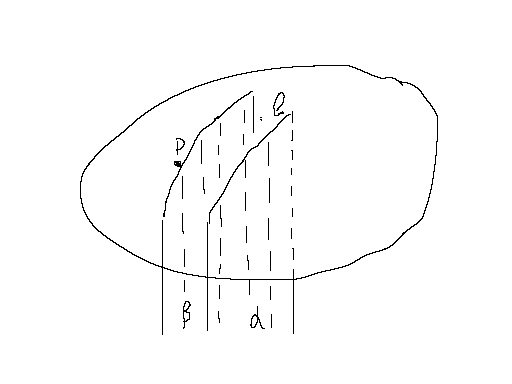}
\begin{center}
Fig.8
\end{center}

\paragraph{Lobachevsky himself uses a slightly different lemma (STP, proposition 28) for the
same purpose . Let me quote only his conclusion (STP, proposition 34), which shows
more precisely the way in which Lobachevsky anticipated the Hilbertean approach:}

\begin{quote}
On the limiting surface [i.e. on the horosphere] sides and angles of triangles hold the
same relations as in the Usual [i.e. Euclidean] geometry
\end{quote}

\paragraph{This crucial observation allowed Lobachevsky to develop (what we call today)
the Hyperbolic trigonometry on the basis of the usual (Euclidean) trigonometry. He writes
down basic principles of this new calculus in the form of four (eqational) identites (STP,
proposition 37, formula 8). In FG and NFG Lobachevsky applies this calculus to a large
class of geometrical problems and in AIG - to calculation of certain integrals, which
earlier were not given any geometrical sense. On the top of that Lobachevsky puts
forward in FG the following general argument purporting to show that the new calculus
guarantees consistency of his Imaginary geometry:}

\begin{quote}
"[1] As far as we are found the equations which represent relations between sides and
angles of triangle ... Geometry turns into Analytics, where calculations are necessarily
coherent and one cannot discover anything what is not already present in the basic
equations. [2] It is then impossible to arrive at contradiction, which would oblige us to
refute first principles, unless this contradiction is hidden in those basic equations
themselves. [3] But one observes that the replacement of sides $a, b, c$ by $ai, bi, ci$
turn these [basic hyperbolic] equations into equations of Spherical Trigonometry. [4]
Since relations between lines in the Usual [i.e. Euclidean] and Spherical geometry are
always the same, [5] the new [i.e. Imaginary] geometry and [Hyperbolic] Trigonometry
will be always in accordance with each other." (FG;  $i$ stands here for the square root
from minus one as usual; numbers in square brackets are introduced for the following analysis of this argument.)
\end{quote}

\paragraph{Let's analyse this complicated argument step by step. First (in part [1]) Lobachevsky claims
that trigonometric relations, which are valid for an arbitrary triangle, allow one to translate the
whole of geometry from the synthetic to the analytic language. He takes this claim for
granted in case of Euclidean geometry and then says that it equally holds in a more
general case of Imaginary geometry. In [2] Lobachevsky assumes that
algebraic transformations are better controllable than synthetic constructive
procedures. The transparency of the ``analytic'' procedures guarantees that if basic
equations contain no hidden contradiction so does the rest. One cannot claim the
same for constructive synthetic procedures since such procedures can bring a
contradiction at any step of reasoning but not only at the initial step of laying out
basic principles. So the analytic means help to reduce the question about consistency
of Imaginary geometry to the same question concerning only foundations of this
geometry, i.e., only the ``basic equation''. [3] is a crucial observation (first made by Lambert) concerning a
profound analogy between Spherical and Hyperbolic geometries. Lobachevsky didn't
understand the sense of this analogy in terms of curvature as we understand it today but took it as a purely formal analogy. His argument, as far as I understand it, is the following. Since spherical
geometry (including spherical trigonometry) is a well-established part of Euclidean
geometry there is no reason to expect any contradiction in it (see [4]). The two
faces of Spherical geometry - the synthetic face and the analytic face - match each other just like in
the case of Plane Euclidean geometry. Hence Spherical trigonometry is consistent. Since
the formal substitution $a\rightarrow ai$ transforms every equation of Spherical trigonometry into an
equation of Hyperbolic trigonometry and since such a formal substitution cannot bring any contradiction into the given theory the Hyperbolic trigonometry is consistent too. But the analytic face of Imaginary geometry (i.e., Hyperbolic trigonometry) matches its synthetic face just like in Spherical case (see [5]). Hence Imaginary geometry (including its synthetic part) is consistent in general. The line of this argument can be pictured with the following diagram:}

\begin{figure}
\includegraphics {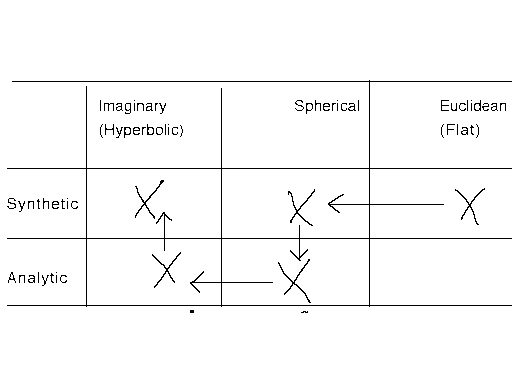}
\end{figure}

\paragraph{It is tempting to interpret this Lobachevsky's argument as a proof of relative consistency in
Hilbert's sense. Even if such reading is not unreasonable one should keep in mind that,
first, this argument is in fact very vaguely formulated and, second, it is produced by
Lobachevsky at the absence of any genuine understanding of what is behind the formal
correspondence between trigonometric identities in Spherical and Hyperbolic cases.
The main source of Lobachevsky's ambiguity here is the lack of any proper distinction
between Imaginary (Hyperbolic), Spherical and Usual (Euclidean) geometries. The
context strongly suggests considering them on equal footing as we do it today. But
remind that Lobachevsky considers the Usual geometry as a special case of
Imaginary and Spherical geometry as a part of Usual. At least the latter assumption is
essential for the argument. When Lobachevsky says that ``relations between lines in
the Usual and Spherical geometry are always the same'' he, in my understanding, looks
at a given sphere as an Euclidean object but not intrinsically. Thus he doesn't think
about it as a model in anything like today's sense. The \emph{formal} character of substitution $a\rightarrow ai$ is obviously explained by Lobachevsky's lack of understanding of its nature. Lobachevsky like Lambert simply noticed the striking analogy between trigonometric identities valid in Imaginary geometry and well-known trigonometric identities for spherical triangles. The analogy suggested considering Plane Imaginary geometry as a sort of Spherical geometry on a sphere of imaginary radius, for example, of radius equal to $i$. Prima facie this didn't make any geometrical sense. But the analogy also suggested that the new trigonometric calculus could work
just as well as Spherical trigonometry, and this in its turn suggested that the
synthetic reasoning behind this new calculus was also correct. This is the core of the
above argument. But obviously the analogy noticed by Lambert and  Lobachevsky was calling for explanation. It is  an irony of history that Lobachevsky's eventual ``formalism'', which was due to the lack of understanding of one particular mathematical question, can be seen today as an anticipation of Hilbert's deliberate formalism based on serious epistmological considerations.}

\paragraph{We have seen that although Lobachevsky had some elements of Hilbertian scheme at
his disposal he was quite strongly attached to the traditional way of geometrical
thinking. Given the historical distance between Lobachevsky and Hilbert this is hardly
surprising. What is more surprising is the unusual way in which Lobachevsky uses
these elements of Hilbertian scheme. The epistemological background behind this scheme suggests no possible reason why Lobachevsky might need a non-standard model of the old good Euclidean plane (rather than some model of Hyperbolic plane) - even if the role of this construction can be easily understood in purely technical terms. This fact shows that the Hilbertian framework is not quite adequate for analyzing Lobachevsky's work. This, once again, is hardly a surprise. But I think that one can get a more interesting lesson from this story by questioning the Hilbertian framework itself. In order to do this let me first introduce some broader historical context.  
}

\section{Extending the context}

 \paragraph{The research focused on the problem of parallels was not the only source of development of Non-Euclidean geometries in 19th century. Another major source was the research of general curve surfaces that first brought significant results in \cite{Gauss:1828}. In this seminal work Gauss first identified geometrical properties of a surface that don't depend on the way, in which this surface is embedded in the 3D Euclidean space. Such properties determine what we call today the \emph{intrinsic geometry} of the given surface. By generalizing Gauss' idea of intrinsic geometry to higher dimensions Riemann in his \cite{Riemann:1854} put forward the notion of (Riemannian) \emph{manifold}. This new generalized notion of geometrical space until today remains on object of active mathematical study.}     
 
 \paragraph{Gauss' long-term interest to the problem of parallels suggests that he already understood that the two lines of research were closely related. However Lobachevsky was apparently wholly unaware about this link. It was first clearly shown by Beltrami in his \cite{Beltrami:1868b} where this author identified Lobachevsky's plane and Lobachevsky's space with Riemannian manifolds of constant negative curvature of dimensions 2 and 3 correspondingly.}

\footnote{In the same year of 1840 when Lobachevsky published his STP Ferdinand Minding
(1806-1885) published in the Crelle Journal a note \cite{Minding:1840}, where he showed
that trigonometrical formulae for triangles formed by geodesics on surfaces of
constant negative curvature can be obtained from trigonometrical formulae for
spherical triangles by replacement of usual trigonometric functions by hyperbolic
functions. Lobachevsky in STP makes a similar observation about straight lines of his
geometry. A communication between the two
authors would most probably lead to the discovery made by Beltrami only in 1868.
However Lobachevsky apparently didn't read this Minding's paper in spite of the fact
that the library of Kazan University had this issue the Crelle Journal. At least the
preserved list of books and journals borrowed by Lobachevsky from the University
Library doesn't include this item \cite{Êèðèìóëëèí:1974}.}

\paragraph{Noticeably this Beltrami's result is less known than another result published by the same author earlier in the same year of 1868. In this earlier paper \cite{Beltrami:1868a} Beltrami claimed that he found a surface (that he called a \emph{pseudosphere}) such that certain curves on this surface (geodesics) behave just like straight lines on the Lobachevsky's plane. To put it in modern (Hilbertain) terms Beltrami claimed that he found an Euclidean model of Lobachevsky's Non-Euclidean geometry. However this result was not quite satisfactory even in Beltrami's own eyes. First, it contained a mistake noticed soon be Hilbert and Helmholz: Beltrami's \emph{pesudosphere} models only a finite region of Lobachevsky's plane but not the whole infinite plane. As Hilbert showed soon this mistake couldn't be corrected. Second and more important, Beltrami realized that the case of Lobachevsky's 3D space cannot be treated similarly. So he looked for a different and more satisfactory solution, which he found and published \cite{Beltrami:1868b} after he read Riemann's memoir \cite{Riemann:1854}. }   

\footnote{From the Hilbertian viewpoint the mistake made by Beltramy in his first paper \cite{Beltrami:1868a} is excusable because this paper also contains what is called today \emph{Beltrami-Klein model} of Lobachevsky's plane.}    

\paragraph{Using Beltrami's results Klein in two papers \cite{Klein:1871} and \cite{Klein:1873} sketched a theoretical framework, which unified most important geometrical theories of the time, including Projective geometry, Affine geometry, Riemannian geometry (the case of constant curvature) and, finally, the ``Imaginary'' Lobachevsky's geometry. In a more developed and systematic form this unificatory account is presented in Klein's posthumously published lectures \cite{Klein:1928}. This work realizes by different means Lobachevsky's and Bolyai's idea of a generalized geometry independent of P5, which comprises Euclidean geometry as a particular case. It made clear why Lobachabsky's and Bolyai's approaches to building such a general theory didn't work as expected: the Euclidean and the Hyperbolic (Lobachevskian) cases taken together don't produce a viable general theory because such a theory leaves out the Elliptic case, i.e., the case of manifolds of constant \emph{positive} curvature. The \emph{absolute geometry} in Bolyai's sense can be compared with a tentative theory of conics that treats parabolas and hyperbolas but doesn't treat circles and ellipses.}  

\paragraph{Thus although the problem of parallels in its traditional setting indeed suggests a possibility of generalization of Euclidean geometry it doesn't point at any viable generalization of this sort. As Weyl \cite{Weyl:1923} describes the situation in 1923:}    

\begin{quote} 
The question of the validity of the ``fifth postulate'', on which historical development started its attack on Euclid, seems to us nowadays to be a somewhat accidental point of departure. The knowledge that was necessary to take us beyond the Euclidean view was, in our opinion, revealed by Riemann.
\end{quote} 

\paragraph{Although a detailed conceptual analysis of the 19th century geometry is out of the scope of this paper I nevertheless mention here Riemann's and Klein's contributions for making it clear that Hilbert's treatise on foundations of geometry \cite{Hilbert:1899} didn't provide even a minimal theoretical background necessary for doing new geometrical research in the beginning of the 20th century. It had a more limited and more specific task: to demonstrate how an older geometrical theory  (Euclidean geometry)  can be reconstructed with a novel Axiomatic Method. Non-Euclidean geometries (in particular, the Lobachevskian geometry) show up here only as means for proving the mutual independence of the geometrical axioms.}     

\paragraph{Hilbert was, of course, well aware about this fact but hoped that his Axiomatic Method could be equally applied in more advanced theoretical contexts. Already the second edition \cite{Hilbert:1903} of Hilbert's \emph{Grundlagen} contained four Supplements treating some contemporary geometrical issues including issues of Rimannian geometry. A more complete Hilbert-style axiomatization of Riemannian geometry has been proposed in 1932 by Veblen and Whitehead \cite{Veblen:1932}. In my view, this result is somewhat controversial. On the one hand, this axiomatization was clearly a success: it allowed for a significant clarification of the whole subject  and was later widely used for educational and research purposes \cite{Scholz:1999}. On the other hand, a comparison of Hilbert's axiomatization of Euclidean geometry with the axiomatization of Riemannian geometry by Veblen and Whitehead immediately shows that Veblen and Whitehead don't really stick to the Hilbertian formalist approach. In the first five Chapters of their tract the authors provide a traditional contentual account of the subject; axioms of the theory are given only in Chapter 6 while in the last Chapter 7 the authors study some general consequences of the axioms. Perhaps this can be explained on the basis of some pragmatic arguments without questioning the Axiomatic Method itself. However the axiomatic account of the two last Chapters also misses essential features of Hilbert's axiomatic account of Euclidean geometry. Unlike Hilbert Veblen and Whitehead don't explicitly discuss \emph{models} of their axiomatic theory, so the very difference between the theory and its models, which is fundamental for Hilbert's approach, is not present in Veblen and Whitehead's account in an explicit form. Thus in spite of the fact that Veblen and Whitehead were indeed inspired by Hilbert's notion of Axiomatic Method, and tried to apply this method in the Riemannian geometry, in order to produce a useful axiomatization of this geometrical discipline they had to compromise against Hilbert's original approach quite severely. } 

\footnote{In Hilbert's account the Analytic geometry is presented as a particular \emph{model} of the (formal) Euclidean geometry, namely as an \emph{arithmetical} model, so Arithmetic is used here as a meta-theory. In Veblen and Whitehead's account analytic devices (and hence Arithmetic) are directly involved into the axioms of their geometrical theory, so in this case Arithmetic is used as a sub-theory rather than a meta-theory. Since Veblen and Whitehead take Arithmetic (as well as relevant set-theoretic notions) for granted  and since the only geometrical primitive term of their theory is that of \emph{point} the only way of producing alternative models, which remains available for the authors, is to assign to the term ``point'' different semantic values. This indeed opens a possibility of putting a Riemannian structure ``on the top of'' mathematical structures of different sorts. However Veblen and Whitehead don't explore this possibility.}

\paragraph{Today we know how to formalize mathematical theories to the extent with comparison to which Hilbert's original approach to foundations of geometry may look completely ``informal''. What remains controversial is the epistemic significance of this procedure. Without going into a broad philosophical discussion about this issue I shall criticize in the next Section some aspects of Hilbert's Axiomatic Method and point to an alternative method of theory-building. We shall see that this alternative method suggested by some recent developments in mathematics has its roots in the 19th century geometry.} 

\section{Rethinking Hilbertian Framework with Lobachevsky}

 \paragraph{Gauss' notion of \emph{intrinsic geometrical property} of a given surface mentioned in the previous Section makes the very notion of geometrical space \emph{relational} in the following sense. Given such a surface one can think of it (i) in the usual way as a two-dimensional \emph{object} living in the Euclidean 3-space and (ii) as a 2-space on its own rights (characterized by the intrinsic properties of the given surface), which is a home for its points, lines, triangles, etc.. Generally, a geometrical \emph{object} can be described in this context as an embedding of the form $s: B \rightarrow C$ where $B$ is a \emph{type} of the given object and $C$ is a \emph{space} where the given object lives. By generalizing again the picture I shall assume $s$ to be a \emph{map} of any sort and not necessarily an embedding. In the given context to be a \emph{space} means to be a target of certain maps; given an object $s: B \rightarrow C$ one calls $C$ a \emph{space} only w.r.t. $s$ (or w.r.t. some other object living in the same space) but not absolutely.}  
 
\paragraph{This way of thinking about geometrical objects and geometrical spaces straightforwardly applies to the traditional Euclidean geometry. Let me distinguish here between two different notions of Euclidean plane: (i) the domain of Euclidean Planimetry and (ii) an object living in the Euclidean 3-space. Notice that one doesn't need this latter notion for doing the Planimetry, it first appears only in the Stereometry! I shall write ``EPLANE'' for the Euclidean plane in the first sense, and write ``eplane'' for the Euclidean plane in the second sense. I shall also write ESPACE and ISPACE to denote the domains of Euclidean and Lobachevskian Stereometry correspondingly. Then an eplane $e$ can be described as a map:}
 
\begin{quote} 
 $e$: EPLANE $\rightarrow$ ESPACE. 
\end{quote} 
 
\paragraph{Notice that there are many different maps of this form (and so we get many different eplanes living in the ESPACE) while the EPLANE and the ESPACE are unique. All eplanes share  these \emph{two} things: first, they share a space, namely ESPACE, and, second, they share a type, namely EPLANE. It is common to say that a given eplane \emph{instantiates} its type EPLANE; on our account the instantiation amounts to mapping this type into a background space like ESPACE.}  

\paragraph{Let me also write ``ISPACE'' for the domain of Lobachevskian Stereometry  (Letter ``I'' stands for ``Imaginary''). Then Lobachevsky's observation about horospheres mentioned in the Section 4 above can be presented as follows: a horosphere $h$ is a map of the form: }

\begin{quote}
$h$: EPLANE $\rightarrow$ ISPACE         
\end{quote}

\paragraph{Thus a horosphere is also an instantiation of the EPLANE but this time obtained with a different background space.} 

\paragraph{A Hilbertian-minded thinker would think of a given eplane and a given horosphere as self-standing entities, construct an isomorphism between them by mapping points of the eplane to appropriate points of the horosphere and vice-versa, and on this basis claim that the eplane and the horosphere represent the same plane Euclidean \emph{structure}. I claim that such an isomorphism is ill-formed. For the eplane and the horosphere are not self-standing entities. The fact that an eplane lives in the Euclidean 3-space is just as essential as the fact that it is an image of the EPLANE. The fact that a horosphere lives in the Hyperbolic 3-space (rather than elsewhere) is equally essential. There is no way to ``take out'' (through an abstraction or otherwise) an eplane from the Euclidean 3-space or to ``take out'' a horosphere from the Hyperbolic space. Any reasonable notion of map (and in particular of isomorphism, i.e. of reversible map) between an eplane and a horosphere \emph{must} take ESPACE and ISPACE into consideration. The isomorphism suggested by the Hilbertian-minded thinker does not meet this requirement leaving the ESPACE and the ISPACE out the consideration as meta-theoretical constructions. This isomorphism is ill-formed because its source and its target are ill-formed. An eplane and a horosphere are not isomorphic; what they share in common is a type, namely the EPLANE, but not an abstract Euclidean structure.}    

\paragraph{A space and a type are not things that one can represent in one's imagination directly. One needs some objects for it. In order to visualize a type one \emph{instantiates} it by an appropriate object (belonging to this given type), i.e. maps the given type into some available space. In order to visualize a space one populates this given space by appropriate objects, i.e. maps available types into this given space. One cannot imagine the EPLANE in the same sense in which one can imagine an eplane but one may have an intuitive grasp on the EPLANE through geometrical objects living in this space: circles, triangles and the like. All such familiar objects are maps from corresponding types into the EPLANE; in particular a circle is a map of the form:}   

\begin{quote}
$c$: CIRCLE $\rightarrow$ EPLANE         
\end{quote}

\paragraph{The ``view on EPLANE from inside'' is multi-faced just like the ``outside view'' obtained by mapping the EPLANE to different spaces like ESPACE or ISPACE. I can see no epistemic reason for thinking that the inside view tells us more about this space than the outside view (or the other round). Nevertheless these two views are clearly different and certainly should not be confused.}

\paragraph{The fact that the EPLANE is populated by objects of different shapes was, of course, well known long before Lobachevsky. However the fact that the EPLANE itself can be seen ``from the outside'' not only as an eplane but also as a horosphere was discovered only by Lobachevsky. This and other similar discoveries (in particular the discovery of a pseudo-sphere by Beltrami) led Hilbert and many of his contemporaries to thinking that one cannot any longer count on intuition in geometry. Given Lobachevsky's notion of horosphere Hilbert would think that since an eplane and a horosphere are two different intuitive representations of the EPLANE and since these intuitive representations look so differently a ``right'' notion of EPLANE cannot be given through intuition at all but must be construed otherwise, namely by formal methods. Unlike Hilbert and his followers I cannot see that our intuition gets something wrong in this situation. An eplane and a horosphere look so differently because they live in different spaces. What is wrong in the traditional geometrical thinking is not the very notion according to which geometrical concepts come with intuitive representations but the idea that a given representation involves just one geometrical concept. In fact any geometrical representation involves at least two things: one, which is represented (a type), and one, which provides a space for this representation. In the traditional Euclidean Planimetry the representation space (the EPLANE) is fixed. Euclidean Stereometry provides extrinsic images of the EPLANE, which all look similarly because these images are all taken against the same fixed background space, namely against the ESPACE. Lobachevsky and other pioneers of non-Euclidean geometry showed that backgrounds (i.e. representation spaces) may vary. Hilbert and other proponents of formal methods didn't take this possibility seriously and looked for a way of doing geometry without any background. Later the proponents of formal methods began to use such methods with a new fixed background of a more ``abstract'' sort, namely with the background universe of \emph{sets}. The alternative approach that I defend here allows for \emph{variation} of backgrounds and moreover makes the very notion of background strictly relational. Hilbert conceives of geometrical space as a ``system of things'' determined by mutual relations between these ``things''. A geometrical space so construed doesn't depend on its eventual relations with any other geometrical space. The proposed approach makes any talk of a single self-standing space senseless. So the sense of being ``relational'' relevant to the suggested alternative approach is stronger than that relevant to the Hilbertian approach.}

\paragraph{Let me now return to the set of questions discussed in the beginning of this paper: Are there different geometries or a single science of geometry? If there exist mutually incompatible geometries how they may coexist? Where do different geometrical spaces live together? The alternative approach just outlined suggests the following answer: all geometrical spaces live in a ``space'' of \emph{maps} between these spaces. The theory of this ``space'' unifies theories of particular geometrical spaces into is a general theory of geometry.}  

\paragraph{Such a notion that involves a class of ``things'', a class of maps between these ``things'' and an operation of \emph{composition} of these maps, which is a subject to several simple rules, is well-known in today's mathematics under the name of \emph{category} \cite{MacLane:1982}. For example all Riemannian manifolds and all maps between these manifolds form the category $R$ of Riemannian manifolds. There are two ways of thinking about $R$. One may have beforehand a notion of Riemannian manifold and a notion of map between such manifolds and then build $R$ by considering all such manifolds and all such maps. Alternatively one may, first, conceive of an ``abstract'' category (without assuming in advance what exactly it consists of) and then, by specifying the ``behavior'' of its maps, make this category look exactly (with an appropriate notion of exactness) like $R$ obtained by the former method. For obvious reason only the latter method is relevant to foundations and can in principle compete with the Hilbertian approach. The reference to the former method is, in fact, not essential. Although it is natural to ask whether on not the new method can reproduce an earlier known result its real purpose is to produce a new notion of geometrical space, which may settle some earlier posed problems. The present discussion concerns such a general notion of space but not a more specific issue of Riemannian manifolds, which I mention here only as a suggestive example.  } 

\paragraph{I warn the reader that the terminology used in this paper differs from one that is commonly used in the Category theory. Recall that I identify a geometrical \emph{object} with a map from a  \emph{type} to a \emph{space}; a space here is a ``meet'' of maps while a type is a ``join'' of maps; a join of maps may well be a meet of some other maps. In the Category theory such types/spaces are called  \emph{objects} (while maps are also called \emph{morphisms}). In order to justify this change of terminology let me stress that a space is not just an object and an object is not just a thing. An object is a thing represented in a background space. When one talks about objects without specifying the background one uses some background space by default. In fact the terms ``type'' and ``space'' are already used in a similar sense in some more special contexts. My types/spaces (i.e. objects in the sense of the usual terminology) are called (or interpreted as) ``types'' in the categorical logic \cite{Bell:1988} and called (or interpreted as) ``spaces'' in any geometrical category. The proposed terminology reflects both these ways of thinking about categories. It reflects a duality between logic and geometry, which becomes apparent when one thinks about an object as a map from some given type to some given space. A map going into the opposite direction turns the given type into a space and the given space into a type. In order to simplify the terminology I propose to call these type/spaces by the common name  of \emph{unit} (rather than  call them ``objects'' as usual). }    

\paragraph{Recall that Hilbert in his \cite{Hilbert:1899} construes the ESPACE as a ``system of things'' comprising primitive objects of several sorts and primitive relations of several sorts that satisfy certain axioms. The alternative approach that I advocate here amounts, roughly, to the replacement of relations by maps. This may look like a merely technical improvement but in fact it implies a significant epistemic shift. The ESPACE construed as a category is not a formal scheme: eplanes, straight lines, points and all its other inhabitants can be conceived of in the usual intuitive way. Yet in a categorical setting the ESPACE has multiple images (by appropriate mappings) in arithmetical and other categories, which in a standard Hilbertain setting are  described as different \emph{models} of ESPACE. The Hilbertian notion of meta-theory turns to be irrelevant to this new context: in a categorical setting the ``models'' ESPACE belong to an extension of the theory of ESPACE but not to a meta-theory distinguished by its epistemic role. I shall not try to suggest such a precise categorical construction of ESPACE in this paper but leave it as a project for the future work.}  

\section{Conclusion}

\paragraph{We have seen that a popular image of the history of geometry of 19th century based on he Hilbertian view onto this discipline, doesn't  quite fit historical evidences. The idea according to which the discovery of Non-Euclidean geometries implied a break with the geometrical intuition in favor of a more abstract mathematical thinking is anachronistic and reflects a later trend, for which Lobachevsky and other pioneers of the new geometry are not responsible.}

\paragraph{Instead of looking at Lobachevsky's way of translating between different theories as an early incomplete grasp of the modern notion of model I suggested a different approach, which takes the notion of such an "incomplete" translation seriously. Then as it so often happens in history apparent shortcomings of an older work become to look as strokes of genius.
}

\paragraph{Lobachevsky's observation that ``[t]he geometry on the limiting sphere [= the horosphere] is exactly the same as on the [Euclidean] plane'' was indeed a great discovery that allowed Lobachevsky to develop an analytic apparatus appropriate for his brand-new ``imaginary geometry''. However unlike Hilbert and his followers Lobachevsky didn't interpret this observation in the sense that the horosphere and the usual Euclidean plane (or more precisely the \emph{eplane}) are different intuitive images of the same abstract entity, which must be ultimately accounted for without using these or other similar images. Actually one can draw a very different epistemic moral from the same geometrical observation: since such different objects as the horosphere and the eplane have the same intrinsic geometry the intrinsic viewpoint is not sufficient for characterizing these objects; this example demonstrates that the intrinsic and the extrinsic viewpoints are equally significant in geometry, so there is no reason to stress one at the expense of the other. The whole Hilbertian idea of construing a geometrical space as a self-standing ``system of things'' (determined by mutual relations between these things) is misleading because geometrical spaces are collective animals (just like points and straight lines on Hilbert's own account!) and any such space is determined by maps from and to some other spaces.}

\paragraph{I don't claim that Lobachevsky held this modern categorical view on geometry that I am advocating here. To the best of my knowledge Lobachevsky never explicitly stressed the controversy between the extrinsic and the intrinsic viewpoints in geometry (although his expression ``the geometry on the limiting sphere'' shows that he was familiar with the notion of intrinsic geometry of a surface). However unlike Hilbert he doesn't say either anything that rules this modern interpretation out.}      

\paragraph{It is true that the Axiomatic Method as designed by Hilbert was a great achievement that allowed, among other things, to cope with Non-Euclidean geometries in a certain way. It is also true that this method proved effective in mathematics of 20th century. However this success doesn't mean that the development of our understanding of basic mathematical issues may or should stop at this point. There emerge new ways of doing mathematics and new ways of thinking about mathematics, which have a bearing on our understanding of history of the subject. Reciprocally, a careful study of history of mathematics may provide useful hints for making further progress in the pure mathematics and its philosophy. The intellectual revolution triggered by Gauss, Lobachevsky, Riemann and other geometers of 19th century still continues and requires further theoretical efforts.}   

\bibliographystyle{plain}
\extrasrussian 
\bibliography{loba} 
\end{document}